\documentclass[12pt, reqno]{tran-l}
\usepackage{ifthen,amsfonts,amsmath,amssymb,epic,eepic}
\usepackage{epsfig,color}

\makeatother

\newtheorem{Thm}{Theorem}
\newtheorem{Cor}{Corollary}
\newtheorem{Lem}{Lemma}

\numberwithin{equation}{section}

\newcommand{\nn}{{\bf{n}}}
\newcommand{\K}{{\text{K}}}

\def\ZZ{{\bold Z}}
\def\RR{{\bold R}}

\def\CC{{\bold C }}

\newcommand{\e}{{\text {e}}}

\newcommand{\cP}{{\mathcal{P}}}

\newcommand{\eqr}[1]{(\ref{#1})}

\begin{document}

\title[Nonproper examples of embedded minimal 
disks in a ball]
{Embedded 
minimal disks:  Proper versus nonproper - global versus local}

\author{Tobias H. Colding}%
\address{Courant Institute of Mathematical Sciences and Princeton University\\
251 Mercer Street\\ New York, NY 10012 and Fine Hall, Washington
Rd., Princeton, NJ 08544-1000}
\author{William P. Minicozzi II}%
\address{Department of Mathematics\\
Johns Hopkins University\\
3400 N. Charles St.\\
Baltimore, MD 21218}
\thanks{The authors were partially supported by NSF Grants DMS
0104453 and DMS 0104187}

\email{colding@cims.nyu.edu, minicozz@jhu.edu}

\begin{abstract}
We construct a sequence of (compact) embedded minimal disks in a ball 
in $\RR^3$ with boundaries in the 
boundary of the ball and  where the curvatures blow up only at the center.   The sequence converges to a 
limit which is not smooth and not proper.  If instead 
the sequence of embedded disks had boundaries in a sequence of balls with radii tending to infinity, 
then we have shown previously that any limit must be smooth and proper.
\end{abstract}

\maketitle

\section{Introduction}

Consider a sequence of (compact) embedded minimal disks
$\Sigma_i\subset B_{R_i}=B_{R_i}(0)\subset \RR^3$ with $\partial
\Sigma_i\subset \partial B_{R_i}$ and
either:\\
(a) $R_i$ equal to a finite constant.\\
(b) $R_i\to\infty$.\\
We will refer to (a) as the {\it local case} and to
(b) as the {\it global case}.  Recall that a surface
$\Sigma \subset \RR^3$ is said to be properly embedded if it is
embedded and the intersection of $\Sigma$ with any compact subset
of $\RR^3$ is compact.  We say that a lamination or foliation is
proper if each leaf is proper.

\begin{figure}[htbp]
    \setlength{\captionindent}{4pt}
    \begin{minipage}[t]{0.5\textwidth}
    \centering\input{loc0.pstex_t}
    \caption{The limit in a ball of a sequence of degenerating helicoids
    is a foliation by parallel planes. This is
    smooth and proper.}\label{f:f1}
    \end{minipage}\begin{minipage}[t]{0.5\textwidth}
\centering\input{loc1.pstex_t}
    \caption{A schematic picture of the limit in Theorem \ref{t:main}
    which is not smooth and
     not proper (the dotted $x_3$-axis is part of the limit).
The limit contains four multi-valued graphs joined at the $x_3$-axis;
$\Sigma_1^+$, $\Sigma_2^+$ above
the plane $x_3=0$ and $\Sigma_1^-$, $\Sigma_2^-$ below the plane.  Each of
the four spirals into the plane.
}\label{f:f2}
\end{minipage}
\end{figure}

We will be interested in the possible limits of  sequences of
minimal disks  $\Sigma_i$ as above where the curvatures blow up,
e.g., $\sup_{B_1 \cap \Sigma_i} |A|^2 \to \infty$ as $i\to
\infty$. In the global case, Theorem $0.1$ in \cite{CM2} gives a
subsequence converging off a Lipschitz curve to a foliation by
parallel planes; cf. fig. \ref{f:f1}.  In particular, the limit is
a (smooth) foliation which is proper.
 We show here in Theorem \ref{t:main} that smoothness and
properness of the limit can fail in the local case; cf. fig. \ref{f:f2}.

\begin{figure}[htbp]
    \setlength{\captionindent}{4pt}
    \begin{minipage}[t]{0.5\textwidth}
\centering\input{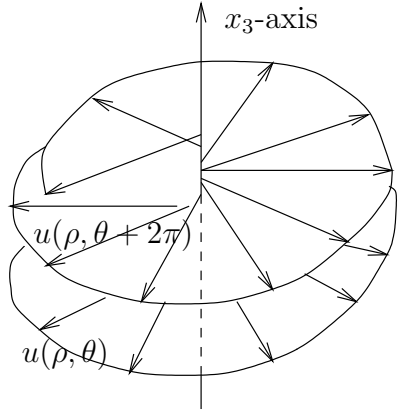}
    \caption{A multi-valued graph of a function $u$.}\label{f:f3}
\end{minipage}
\end{figure}

We will need the notion of a multi-valued graph; see fig. \ref{f:f3}. Let
$D_r \subset \CC$ be the disk in the plane centered at the origin
and of radius $r$ and let $\cP$ be the universal cover of the
punctured plane $\CC\setminus \{0\}$ with global polar coordinates
$(\rho, \theta)$ so $\rho>0$ and $\theta\in \RR$.  An {\it
$N$-valued graph} on the annulus $D_s\setminus D_r$ is a single
valued graph of a function $u$ over $\{(\rho,\theta)\,|\,r<
\rho\leq s\, ,\, |\theta|\leq N\,\pi\}$.

In Theorem \ref{t:main}, we  construct a sequence of disks
$\Sigma_i \subset B_1 = B_1 (0)\subset \RR^3$ as above where
 the curvatures blow up only at $0$ (see 
(1) and 
(2)) and
$\Sigma_i \setminus \{ {\text{$x_3$-axis}} \}$ consists of two
multi-valued graphs for each $i$; see  (3).
Furthermore (see (4)),
$\Sigma_i  \setminus \{ x_3 = 0 \}$  converges to
two
  embedded minimal disks $\Sigma^- \subset \{ x_3 < 0 \}$
and  $\Sigma^+ \subset \{ x_3 > 0 \}$ each of which spirals
 into $\{ x_3 =
0 \}$ and thus is not proper; see fig. \ref{f:f2}.

\begin{Thm} \label{t:main}
There is  a sequence  of compact embedded minimal disks $0 \in
\Sigma_i \subset  B_1 \subset \RR^3$
 with $\partial \Sigma_i \subset
\partial B_1$ and
containing the vertical segment
$\{ (0,0,t) \, | \, |t|<1 \} \subset \Sigma_i$ so:
\begin{enumerate}
\item[(1)] $\lim_{i\to \infty} |A_{\Sigma_i}|^2 (0) = \infty$. 
\label{i:1}

\item[(2)]
$\sup_i \sup_{\Sigma_i \setminus B_{\delta} } |A_{\Sigma_i}|^2 < \infty$
for all $\delta > 0$.
\label{i:2}

\item[(3)]
$\Sigma_i \setminus \{ {\text{$x_3$-axis}} \} =
\Sigma_{1,i} \cup \Sigma_{2,i}$ for multi-valued graphs
$\Sigma_{1,i}$ and
$\Sigma_{2,i}$.
\label{i:3}

\item[(4)] 
$\Sigma_i  \setminus \{ x_3 = 0 \}$  converges to     two
embedded minimal disks $\Sigma^{\pm} \subset \{ \pm x_3 > 0 \}$
with $\overline{\Sigma^{\pm}} \setminus \Sigma^{\pm} = B_1 \cap \{
x_3 = 0\}$. Moreover, $\Sigma^{\pm} \setminus \{
{\text{$x_3$-axis}} \} = \Sigma_1^{\pm} \cup \Sigma_2^{\pm}$ for
multi-valued graphs $\Sigma_{1}^{\pm}$ and $\Sigma_{2}^{\pm}$ each of
which spirals into $\{ x_3 =
0 \}$; see fig. \ref{f:f2}.
\label{i:4}
\end{enumerate}
\end{Thm}

It follows from (4)
that $\Sigma_i \setminus \{ 0 \}$ converges to a lamination
of $B_1 \setminus \{ 0 \}$    (with leaves $\Sigma^-$, $\Sigma^+$,
and $B_1 \cap \{ x_3 = 0 \} \setminus \{ 0 \}$) which
 does not extend to a lamination of $B_1$.
Namely, $0$ is not a removable singularity.

  The multi-valued graphs that
we will consider will never close up; in fact they will all be
embedded.
The most important example of an embedded minimal multi-valued
graph comes from the helicoid.  The {\it helicoid} is the minimal
surface $\Sigma$ in $\RR^3$ parametrized by $(s\cos t,s\sin t,t)$
where $s,\,t\in \RR$. Thus $\Sigma\setminus
\{ {\text{$x_3$-axis}} \}=\Sigma_1\cup \Sigma_2$,
 where $\Sigma_1$, $\Sigma_2$ are  $\infty$-valued graphs on
$\CC\setminus \{0\}$. $\Sigma_1$ is the graph of the function
$u_1(\rho,\theta)=\theta$ and $\Sigma_2$ is the graph of the
function $u_2(\rho,\theta)=\theta+\pi$.

\vskip2mm
 We will use standard $(x_1 , x_2 , x_3)$ coordinates on
$\RR^3$ and $z=x+i\, y$ on $\CC$.  Given $f:\CC \to \CC^n$,
$\partial_x f$ and $\partial_y f$ denote $\frac{\partial
f}{\partial x}$ and $\frac{\partial f}{\partial y}$, respectively;
similarly, $\partial_z f =  (\partial_x f - i \partial_y f) /2$.
 For $p \in \RR^3$ and $s
> 0$, the
 ball in $\RR^3$ is
$B_s(p)$. $\K_{\Sigma}$ is the sectional curvature of a smooth
surface $\Sigma$. When $\Sigma$ is immersed in $\RR^3$, then
$A_{\Sigma}$ will
be its second fundamental form (so when $\Sigma$ is minimal, then
$|A_{\Sigma}|^2=-2\,\K_{\Sigma}$). When $\Sigma$ is oriented,
$\nn_{\Sigma}$ is the unit normal.

\section{Preliminaries on the Weierstrass representation}

Let $\Omega \subset \CC$ be a domain. The  classical Weierstrass
representation (see \cite{Os}) starts from a meromorphic function
$g$ on $\Omega$  and a holomorphic one-form $\phi$ on  $\Omega$
and associates a (branched) conformal minimal immersion $F: \Omega
\to \RR^3$ by
\begin{equation}    \label{e:ws1}
    F(z) = {\text{Re }} \int_{\zeta \in \gamma_{z_0,z}}
\left( \frac{1}{2} \, (g^{-1} (\zeta) - g (\zeta) )
    , \frac{i}{2} \, (g^{-1} (\zeta)
    +g (\zeta) ) , 1 \right) \, \phi (\zeta) \, .
\end{equation}
Here $z_0 \in \Omega$ is a fixed base point and the integration is
along a path $\gamma_{z_0,z}$ from $z_0$ to $z$.
The choice of $z_0$
changes $F$ by adding a constant. We will assume that
$F(z)$ does not depend on the
choice of path $\gamma_{z_0,z}$; this is the case,
for example, when $g$ has no zeros or poles
and $\Omega$ is simply connected.

 The unit normal $\nn$
and Gauss curvature $\K$ of the resulting surface are then (see
sections $8$, $9$ in \cite{Os})
\begin{align}    \label{e:nn}
    \nn  &= \left( 2 \, {\text{Re }} g , 2 \, {\text{Im }} g , |g|^2 -1
    \right)/ (|g|^2 +1)
    \, , \\
    \K &= - \left[ \frac{4 |\partial_z g| \, |g|}{|\phi| \, (1 +|g|^2)^2}
    \right]^2
     \, . \label{e:K}
\end{align}
Since the pullback $F^{\ast} (dx_3)$ is ${\text{Re }} \phi$ by \eqr{e:ws1},
$\phi$ is usually called the {\it height differential}.
By \eqr{e:nn},  $g$ is the composition
of the Gauss map followed by stereographic projection.

To ensure that $F$ is an immersion (i.e., $dF \ne 0$),
we will assume that
 $\phi$ does not vanish and $g$ has no zeros or poles. The
two standard examples are
 \begin{align}
 &g (z) = z, \, \phi (z)  = dz/z ,
\, \Omega = \CC \setminus \{ 0 \} {\text{ giving a catenoid}}
\, , \\
  &g (z) = \e^{iz} , \,
\phi (z) = dz , \, \Omega = \CC  {\text{ giving a helicoid}} \, .
\label{e:hel}
\end{align}

The next lemma records the differential of  $F$.

\begin{Lem} \label{l:gammax}
If $F$ is given by \eqr{e:ws1} with $g (z) = \e^{i \, (u(z)+i
v(z))}$ and $\phi= dz$, then
\begin{align}
 \partial_x F
    &=  ( \sinh v  \, \cos u ,  \sinh v  \, \sin u  , 1)  \label{e:ws6x} \,
    , \\
    \partial_y F
    &= (\cosh v \,  \sin u , - \cosh v \,  \cos u , 0)  \, . \label{e:ws6}
\end{align}
\end{Lem}

\section{The proof of Theorem \ref{t:main}}

To show Theorem \ref{t:main},
we first construct a one-parameter  family (with parameter $a\in
(0,1/2)$) of minimal immersions $F_a$ by making a specific choice
of  Weierstrass data $g = \e^{i h_a}$ (where $h_a = u_a+i\, v_a$),
$\phi = dz$, and domain $\Omega_a$ to use in \eqr{e:ws1}.
 We show in Lemma \ref{l:dg1a} that
this one-parameter family of immersions is compact.
Lemma \ref{l:dg1}
shows
that the immersions $F_a : \Omega_a \to \RR^3$ are embeddings.

\begin{figure}[htbp]
    \setlength{\captionindent}{4pt}
    \begin{minipage}[t]{0.5\textwidth}
    \centering\input{loc3.pstex_t}
    \caption{The domain $\Omega_a$.}\label{f:f4}
    \end{minipage}\begin{minipage}[t]{0.5\textwidth}
    \centering\input{loc3a.pstex_t}
    \caption{$\Omega_0=\cap_{a>0}\Omega_a \setminus \{ 0\}$ and
    its two components $\Omega^+_0$ and $\Omega^-_0$.}\label{f:f5}
    \end{minipage}
\end{figure}

  For each $0<a<1/2$,  set (see fig. \ref{f:f4})
\begin{equation}   \label{e:dg0}
    h_a (z)  = \frac{1}{a} \,
\arctan \left( \frac{z}{a} \right) {\hbox{ on }}
    \Omega_a = \{ (x,y) \, | \, |x| \leq 1/2 ,
\, |y| \leq (x^2 + a^2)^{3/4}/2  \}
 \, .
\end{equation}
Note that $h_a$ is well-defined since $\Omega_a$ is simply
connected and $\pm i \, a \notin \Omega_a$.
 For future reference
\begin{align}   \label{e:dg1}
    \partial_z h_a (z)&= \frac{1}{z^2 + a^2}
=  \frac{x^2 +a^2 - y^2 -2i \, xy}{(x^2 +a^2 - y^2)^2
    + 4x^2y^2} \, , \\
    \K_a (z) &= \frac{-|\partial_z h_a|^2}{\cosh^4 v_a}
    =  \frac{-|z^2 + a^2|^{-2}}{\cosh^{4} \left( {\text{Im }}
    \arctan ( z/a) / a \right)}  \, .  \label{e:K2}
\end{align}
Here
\eqr{e:K2} used \eqr{e:K}.
Note that, by the Cauchy-Riemann equations,
\begin{equation}  \label{e:creqns}
    \partial_z h_a =  (\partial_x - i \, \partial_y) (u_a + i v_a)/2 =
 \partial_x u_a - i \,\partial_y u_a = \partial_y v_a + i \,\partial_x v_a \, .
\end{equation}

In the rest of this paper we let
$F_a : \Omega_a \to \RR^3$ be from \eqr{e:ws1} with $g =
\e^{i \, h_a}$,  $\phi = dz$, and  $z_0 =0$.
Set $\Omega_0 = \cap_a \Omega_a \setminus \{ 0\}$, so $\Omega_0 =
\{ (x,y) \, |\, 0<|x| \leq 1/2 , \, |y| \leq |x|^{3/2} /2 \}$;
see fig. \ref{f:f5}. The
family of functions $h_a$ is not compact since $\lim_{a\to 0}
|h_a| (z) = \infty$ for $z\in \Omega_0$.  However, the next lemma
shows that the family of immersions $F_a$ is compact.

\begin{Lem} \label{l:dg1a}
If $a_j \to 0$, then there is a
subsequence $a_i$ so $F_{a_i}$  converges uniformly in $C^2$ on
compact subsets of $\Omega_0$.
\end{Lem}

\begin{proof}
Since
$h_a$ and $-1/z$ are holomorphic and
\begin{equation}  \label{e:dcv}
    \left| \partial_z h_a (z) - \partial_z (-1/z) \right|
=  a^2 \, |z|^{-2} \, |z^2 +a^2|^{-1} \, ,
\end{equation}
we get easily that $\nabla h_a$ converges
as $a\to 0$ to $\nabla (-1/z)$
uniformly on compact subsets of $\Omega_0$.
Since each $v_a (x,0)=0$, the fundamental theorem of
calculus gives that the $v_a$'s converge uniformly in $C^1$
on compact subsets of $\Omega_0$.
(Unfortunately, the $u_{a}$'s do not converge.)

Let $\Omega_0^{\pm} = \{ \pm x > 0 \} \cap \Omega_0$  be the two
components of $\Omega_0$; see fig. \ref{f:f5}.
Set $b_j^{+} = u_{a_j} (1/2)$ and
$b_j^{-} = u_{a_j} (- 1/2)$
and
choose a subsequence $a_i$ so both $b_i^{-}$ and $b_i^{+}$
converge modulo $2 \pi$ (this is possible since $T^2 = \RR^2 /
(2\pi \ZZ^2)$ is compact). Arguing as above,  $h_{a_i} -
b_i^{\pm}$ converges uniformly in $C^1$ on compact subsets of
$\Omega_0^{\pm}$. Therefore, by Lemma \ref{l:gammax}, the minimal
immersions corresponding to Weierstrass data $g = \e^{i (h_{a_i} -
b_i^{\pm})}$, $\phi = dz$ converge uniformly in $C^2$ on compact
subsets of $\Omega_0^{\pm}$ as $i\to \infty$.
\end{proof}

\begin{figure}[htbp]
    \setlength{\captionindent}{4pt}
    \begin{minipage}[t]{0.5\textwidth}
\centering\input{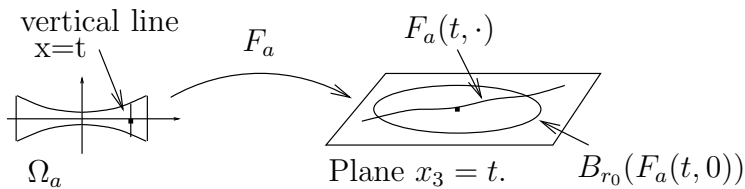}
    \caption{A horizontal slice in Lemma \ref{l:dg1}.}\label{f:f6}
\end{minipage}
\end{figure}

\begin{figure}[htbp]
    \setlength{\captionindent}{4pt}
    \begin{minipage}[t]{0.5\textwidth}
\centering\input{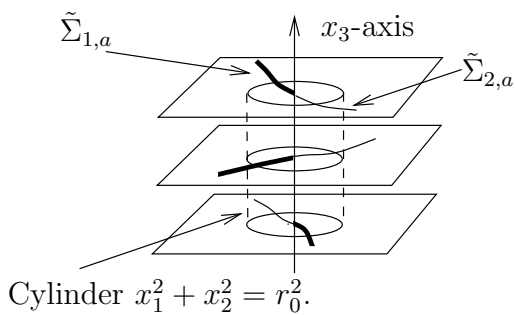}
    \caption{Horizontal slices of $F_a (\Omega_a)$ in Lemma \ref{l:dg1}.} \label{f:f7}
\end{minipage}
\end{figure}

The main difficulty in proving Theorem \ref{t:main} is showing
that the immersions $F_a : \Omega_a \to \RR^3$ are embeddings.
This will follow easily from (A) and (B) below.  Namely, we show
in  Lemma \ref{l:dg1},  see fig.
\ref{f:f6} and \ref{f:f7}, that for $|t| \leq 1/2$:

\noindent (A) The horizontal slice $\{x_3 = t \} \cap F_a
(\Omega_a)$ is the image of the vertical segment $\{ x = t \}$ in
the plane, i.e., $x_3 (F_a(x,y)) = x$;
 see
\eqr{e:dg2a}.

\noindent (B) The image  $F_a \left( \{ x = t \} \cap \Omega_a
\right)$ is a graph over a line segment in the plane $\{ x_3 = t
\}$ (the line segment will depend on $t$); see \eqr{e:dg2}.

\noindent (C) The boundary of the  graph in (B)
 is outside the ball $B_{r_0} ( F_a
( t , 0 ))$  for some $r_0 > 0$ and all $a$;  see \eqr{e:dg3}.

\begin{Lem} \label{l:dg1}
\begin{align}  \label{e:dg2a}
    &x_3 ( F_a (x,y) ) = x  \, . \\
 \label{e:dg2}
    &\hbox{The curve } F_a(x,\cdot) : [-(x^2 + a^2)^{3/4}/2 ,
(x^2 + a^2)^{3/4}/2] \to \{ x_3 = x \}
    \hbox{ is a graph} \, . \\
     &|F_a(x , \pm (x^2 + a^2)^{3/4}/2) - F_a(x,0)| > r_0\text{ for some } r_0>0 \text{ and all } a \, .
\label{e:dg3}
\end{align}
\end{Lem}

\begin{proof}
Since $z_0 = 0$ and $\phi = dz$, we get \eqr{e:dg2a} from \eqr{e:ws1}.  Using
$y^2 < (x^2 + a^2)/4$ on $\Omega_a$,
 \eqr{e:dg1} and \eqr{e:creqns} give
\begin{align}   \label{e:dg4}
    | \partial_y u_a (x,y) | &=   \frac{2\, |xy|}{(x^2 +a^2 - y^2)^2
    + 4x^2y^2} \leq  \frac{4\, |xy|}{(x^2 +a^2)^2}\, , \\
    \partial_y v_a (x,y) &= \frac{x^2 +a^2 - y^2}{(x^2 +a^2 - y^2)^2
    + 4x^2y^2} > \frac{3}{8(x^2 +a^2)}
    \, . \label{e:dg5}
\end{align}
Set $y_{x,a} = (x^2 + a^2)^{3/4}/2$.
Integrating \eqr{e:dg4} gives
\begin{equation}    \label{e:dg7}
    \max_{|y| \leq  y_{x,a}} |u_a (x,y) - u_a (x,0)| \leq
\int_0^{ y_{x,a} } \frac{4|x|t}{(x^2 +a^2)^2} \, dt
    = \frac{|x|}{ 2\,(x^2 +a^2)^{1/2} }
    < 1 \, .
\end{equation}
Set $\gamma_{x,a} (y) = F_a(x,y)$.
Since $v_a(x,0) = 0$ and $\cos (1) >
1/2$, combining
\eqr{e:ws6} and \eqr{e:dg7} gives
\begin{equation}    \label{e:dg8}
     \langle \gamma_{x,a}' (y) , \gamma_{x,a}'(0) \rangle =
    \cosh v_a(x,y) \,  \cos (u_a(x,y) - u_a(x,0) )
    > \cosh v_a(x,y) / 2  \, .
\end{equation}
Here $\gamma_{x,a}' (y) = \partial_y F_{a} (x,y)$.
By \eqr{e:dg8},  the angle between
$\gamma_{x,a}' (y)$ and $\gamma_{x,a}'(0)$ is always less than $\pi /2$;
this gives \eqr{e:dg2}.  Since
$v_a(x,0)=0$, integrating \eqr{e:dg5} gives
\begin{equation}    \label{e:dg9}
    \min_{ y_{x,a}/2 \leq   |y| \leq y_{x,a} } |v_a(x,y)| \geq
\int_0^{ \frac{y_{x,a}}{2} }
   \frac{3 \, dt}{8(x^2 +a^2)}
    = \frac{3}{32 (x^2 + a^2)^{1/4}}  \, .
\end{equation}
Integrating \eqr{e:dg8} and using \eqr{e:dg9} gives
\begin{equation}    \label{e:dg10}
     \langle \gamma_{x,a} (y_{x,a})
    - \gamma_{x,a}(0) , \gamma_{x,a}'(0) \rangle >
    \frac{(x^2 + a^2)^{3/4}}{16} \, \e^{  (x^2 + a^2)^{-1/4}/11  }
     \, .
\end{equation}
Since $\lim_{s\to 0} s^{3} \, \e^{  s^{-1}/11 } = \infty$, \eqr{e:dg10}
and its analog for $\gamma_{x,a} (-y_{x,a})$ give
\eqr{e:dg3}.
\end{proof}

\begin{Cor} \label{c:dg1}
See fig. \ref{f:f7}.  Let $r_0$ be given by \eqr{e:dg3}.
\begin{enumerate}
 \item[(i)] 
$F_a$ is an embedding.
 \item[(ii)] 
 $F_a (t,0) =
(0,0,t)$ for $|t| < 1/2$. 
 \item[(iii)] 
$\{ 0< x_1^2 +
x_2^2 < r_0^2 \} \cap F_a (\Omega_a) = \tilde{\Sigma}_{1,a} \cup
\tilde{\Sigma}_{2,a}$
 for multi-valued graphs
   $ \tilde{\Sigma}_{1,a} , \, \tilde{\Sigma}_{2,a}$ over $D_{r_0} \setminus \{ 0 \}.$
\end{enumerate}
\end{Cor}

\begin{proof}
Equations \eqr{e:dg2a} and \eqr{e:dg2} immediately give (i).

Since  $z_0 = 0$, $F(0,0)=(0,0,0)$.  Integrating \eqr{e:ws6x}
and using $v_a(x,0)=0$ then gives (ii).

 By
\eqr{e:nn}, $F_a$ is ``vertical,'' i.e., $\langle \nn ,
(0,0,1) \rangle = 0$, when $|g_a| = 1$.  However, $|g_a(x,y)|=1$
exactly when $y=0$ so that, by (ii), the image is graphical away
from the $x_3$-axis.
 Combining this with \eqr{e:dg3} gives (iii).
\end{proof}

Corollary \ref{c:dg1} constructs the embeddings $F_a$ that will be
used in Theorem \ref{t:main} and shows property (3).
To prove
Theorem \ref{t:main}, we need therefore only
show  (1), (2), and (4).

\begin{proof}
(of Theorem \ref{t:main}). By scaling, it suffices to  find a
sequence $\Sigma_i \subset B_{R}$ for some $R > 0$.
 Corollary  \ref{c:dg1} gives minimal embeddings $F_a: \Omega_a \to \RR^3$
with $F_a (t,0) = (0,0,t)$
for $|t| < 1/2$ and so (3) 
holds for any $R \leq r_0$.
Set $R = \min \{ r_0/2 , 1/4 \}$ and $\Sigma_{i} = B_{R} \cap F_{a_i}
(\Omega_{a_i})$,  where the sequence $a_i$ is to be determined.

To get (1),
simply note that, by \eqr{e:K2},
$|\K_a| (0)  =   a^{-4} \to \infty$ as $a\to 0$.

We next show (2).
First, by \eqr{e:K2}, $\sup_a \sup_{ \{ |x| \geq
\delta \} \cap \Omega_a } |\K_a|  < \infty$ for all $\delta > 0$.
Combined  with (3) and Heinz's curvature estimate for minimal
graphs (i.e.,  11.7 in \cite{Os}), this gives (2).

To get (4),  use Lemma \ref{l:dg1a} to choose $a_i \to 0$ so the
mappings $F_{a_i}$ converge uniformly in $C^2$ on compact subsets
to $F_0 : \Omega_0 \to \RR^3$.  Hence, by Lemma \ref{l:dg1},
$\Sigma_i \setminus \{x_3 = 0 \}$ converges to two embedded
minimal disks $\Sigma^{\pm} \subset F_0 (\Omega_0^{ \pm})$ with
$\Sigma^{\pm} \setminus \{ {\text{$x_3$-axis}} \} = \Sigma^{\pm}_1
\cup \Sigma^{\pm}_2$ for multi-valued graphs $\Sigma^{\pm}_j$. To
complete the proof, we  show that each graph $\Sigma^{\pm}_j$ is
$\infty$-valued (and, hence, spirals into $\{ x_3 =0 \}$). Note
that, by (3) and \eqr{e:ws6}, the level sets $\{ x_3 = x \} \cap
\Sigma^{\pm}_j$ are graphs over the line in the direction
\begin{equation}    \label{e:infv0}
    \lim_{a \to 0} (\sin u_a (x,0) , - \cos u_a (x,0) , 0) \, .
\end{equation}
Therefore, since an easy calculation gives for $0< t< 1/4$ that
\begin{equation}    \label{e:infv}
    \lim_{a\to 0} |u_a (t,0) - u_a (2t,0)| = 1/(2t)
    \, ,
\end{equation}
we see that $\{t < |x_3| < 2t \} \cap \Sigma^{\pm}_j$ contains an
embedded $N_t$-valued graph where $N_t \approx 1/(4\pi t) \to
\infty$ as $t\to 0$ . It follows that $\Sigma^{\pm}_j$
 must  spiral into
$\{ x_3 =0 \}$, completing (4).
\end{proof}


\begin{thebibliography}{999}

\bibitem[CM1]{CM1}
T.H. Colding and W.P. Minicozzi II,
Embedded minimal disks, To appear in The Proceedings of the
Clay Mathematics Institute Summer School on the Global
Theory of Minimal Surfaces.  MSRI. math.DG/0206146.
\bibitem[CM2]{CM2}
\bysame, The space of embedded minimal
surfaces of fixed genus in a $3$-manifold IV; Locally simply connected,
preprint, math.AP/0210119.
\bibitem[Os]{Os}
R. Osserman,
A survey of minimal surfaces,
{\it Dover},  2nd. edition (1986).

\end{thebibliography}
\end{document}